\documentclass{article}
\usepackage[cp1251]{inputenc}
\usepackage[russianb]{babel}
\usepackage{latexsym,amscd,amssymb,amsopn,amsthm,amsfonts,amsmath}

\newtheorem{theorem}{Теорема}
\newtheorem{corollary}{Следствие}
\newtheorem{lemma}{Лемма}
\newtheorem{notation}{Замечание}

\DeclareMathOperator{\Int}{Int} 

\begin{document}

\title{АСИМПТОТИКА КОЭФФИЦИЕНТОВ ТЕЙЛОРА РАЦИОНАЛЬНОЙ ФУНКЦИИ ДВУХ ПЕРЕМЕННЫХ С ЛИНЕЙНЫМ
МНОЖЕСТВОМ ОСОБЕННОСТЕЙ}
\author{A.П. Ляпин\footnote{\copyright\ А.П. Ляпин, Красноярский государственный
университет, 2006, lalex@krasu.ru}}
\date{}
\maketitle

Пусть $x=(x_1,\ldots,x_n)$ --- точка $n$-мерной целочисленной
решетки $\mathbb{Z}^n$ и $\mathbb{Z}^n_+$ --- подмножество этой
решетки, состоящее из точек с целочисленными неотрицательными
координатами. Обозначим $z=(z_1, \ldots,z_n)$ точку $n$-мерного
комплексного пространства $\mathbb{C}^n$.

Производящая функция $F(z)=\sum_{x\in \mathbb{Z}^n_+}f(x)z^x$
служит эффективным средством изучения последовательности
$\{f(x)\}, x\in\mathbb{Z}^n_+$ в различных задачах математического
анализа и комбинаторики. При изучении ее асимптотического
поведения простейшим будет случай, когда $F(z)$ --- рациональная
функция, причем для $n=1$ из основной теоремы о вычетах следует,
что асимптотика $f(x)$ определяется ближайшим полюсом функции
$F(z)$. В случае многомерной последовательности возникает ряд
трудностей принципиального характера, одной из которых считается
отсутствие понятия кратного асимптотического ряда. Один из
способов исследования асимптотического поведения кратной
последовательности $f(x)$,  $x=(x_1,\ldots,x_n) \in
\mathbb{Z}^n_+$
--- это изучение асимптотики на <<диагоналях>> $x=kp$, где $p$ ---
фиксированная точка из $\mathbb{Z}^n_+$, а $k=0,1,2,\ldots$. Такой
подход применялся в \cite{tz91} в связи с решением проблемы
устойчивости двумерных цифровых рекурсивных фильтров, имеющих
рациональную передаточную функцию. При этом был предложен метод
получения асимптотических оценок коэффициентов Тейлора в случае,
когда множество особенностей рациональной функции двух переменных
гладкое и пересекается с остовом соответствующего полицилиндра в
конечном числе точек.

В работах \cite{orlov2} и \cite{orlovN} методы работы \cite{tz91}
применялись для получения асимптотических оценок коэффициентов
Тейлора алгебраических и мероморфных функций двух переменных. В
случае функций с полюсами на объединении конечного числа
гиперплоскостей  оценки на коэффициенты Тейлора типа
<<$O$-большое>> впервые приведены для $n=2$ в \cite{makosi}, а для
$n>2$ рассмотрены в \cite{orlovN}. Этот случай представляет
большой интерес с точки зрения перечислительного комбинаторного
анализа. Например, в задачах о числе решеточных путей и задаче о
покрытии костями домино соответствующая производящая функция
представляет собой рациональную функцию двух переменных,
знаменатель которой --- произведение линейных множителей
(\cite{pem}).

В данной работе при некоторых ограничениях на линейные множители
знаменателя рациональной функции двух переменных найден главный
член асимптотического разложения коэффициентов $f(kp),
k\to\infty$.


Пусть $F(z,w)$ --- рациональная функция двух комплексных
переменных $(z,w)\in \mathbb{C}^2$, регулярная в начале координат,
и $f(x,y)$
--- коэффициенты ее разложения в ряд Тейлора
\begin{align} \label{f1}
F(z,w)=\frac{P(z,w)}{Q(z,w)}=\sum_{x,y\geqslant 0} f(x,y)z^xw^y,
\end{align}
причем $Q(0,0)\neq 0$ и дробь $\frac{P}{Q}$ --- несократимая.
Рассмотрим ситуацию при следующих условиях:
\begin{itemize}
    \item[(\textbf{I})]знаменатель функции $F$ разлагается в произведение линейных множителей
        \begin{align*}
            Q=\prod_{i=1}^{m}Q_i=\prod_{i=1}^{m}(1-a_iz-b_iw),
        \end{align*}
        т.е. множество особых точек $\mathcal{V}=\{Q(z,w)=0\}$ функции \eqref{f1} является объединением
        системы комплексных прямых  $\mathcal{V}=\cup \mathcal{V}_i$, где $\mathcal{V}_i=\{Q_i=0\}, i=1,\ldots,m$;
    \item[(\textbf{II})]система прямых $\{\mathcal{V}_i\}_{i=1}^{m}$ находится в общем положении в том смысле, что
        любые две из них пересекаются, а любые три - не пересекаются;
    \item[(\textbf{III})]коэффициенты $a_i,b_i$ многочленов $Q_i,\, i=1,\ldots,m$ --- положительные числа.
\end{itemize}

Обозначим $Q_0=z$ и $Q_{m+1}=w$. Пусть многоугольник $M$ определен
в положительном октанте $\mathbb{R}^2_+$ системой линейных
неравенств $Q_i(z,w)>0,\, i=0,\ldots,m+1$. Из условия
\textbf{(III)} следует, что $M$ не пусто. Пусть $(z_i,w_i)$ ---
решение системы уравнений $Q_i(z,w)=qzQ_z(z,w)-pwQ_w(z,w)=0$ для
$i=1,\ldots,m$, а $(z_{ij},w_{ij})$  --- решение системы линейных
уравнений $Q_i(z,w)=Q_j(z,w)=0$ с определителем
$\Delta_{ij}=a_ib_j-a_jb_i \neq 0, \,i\neq j$, и $i,j =
1,\ldots,m$.

Определим в положительном октанте $\mathbb{R}^2_+$ множества двух
видов:
\begin{equation}\begin{split} \label{cones}
K_i&=\{(p,q)\in \mathbb{R}^{2}_{+}:\underset{M}{\max}\,
{z^pw^q}=\underset{\mathcal{V}_i}{\max}\,z^pw^q=z_i^pw_i^q\}, \\
\Omega_{ij}&=\{(p,q)\in
\mathbb{R}^{2}_{+}:\underset{M}{\max}\,z^pw^q=z_{ij}^pw_{ij}^q\}.
\end{split}\end{equation}
Множества $K_i, \Omega_{ij}$ представляют собой конусы,
объединение которых совпадает с $\mathbb{R}^2_+$, при этом
некоторые их них могут оказаться пустыми, а непустые конусы могут
пересекаться лишь по граничным лучам. Можно утверждать, что почти
всякая точка $(p,q)\in \mathbb{Z}^2_{+}$ принадлежит внутренности
одного из <<целочис-ленных>> конусов вида $\mathbb{Z}^2_{+}\cap
K_i$ и $\mathbb{Z}^2_{+}\cap \Omega_{ij}$.

Множество $M$ в общем случае представляет многоугольник в
$\mathbb{R}^2_+$, две стороны которого\, лежат на координатных
осях, а остальные --- на некоторых (или всех) прямых
$\mathcal{V}_i$. Не умаляя общности, этими прямыми можно считать
$\mathcal{V}_1, \ldots, \mathcal{V}_r,\, r\leqslant m$. Вершины
многоугольника $M$, не совпадающие с точкой (0,0), имеют
координаты $(z_{i,i+1},w_{i,i+1}),\, i=0,...,r$, при этом будут
выполняться неравенства
$\frac{a_1}{b_1}<\frac{a_2}{b_2}<\frac{a_3}{b_3}<...<\frac{a_r}{b_r}$.
Построим систему векторов $\eta_{10} = (a_1z_{01}, b_1w_{01}) =
(1, 0),\, \eta_{11} = (a_1z_{12}, b_1w_{12}),\, \eta_{21} =
(a_2z_{12}, b_2w_{12}),\, \eta_{22} = (a_2z_{23},\, b_2w_{23}),
...,\, \eta_{n,n-1} = (a_nz_{n-1,n}, b_nw_{n-1,n}),\, \eta_{n,n} =
(a_nz_{n,n+1}, b_nw_{n,n+1}) = (0,1)$, тогда конусы двух видов,
построенные на этих векторах, совпадут с определенными формулами
\eqref{cones}, а их представление через образующие вектора
$\eta_{ij}$ будет иметь вид
\begin{align*}
K_i &= \{(p,q):(p,q)=\alpha\eta_{i,i-1}+\beta\eta_{i,i},\, \alpha,\beta > 0\},\,i=1, \ldots,n,\\
\Omega_{i, i+1} &=
\{(p,q):(p,q)=\alpha\eta_{i,i}+\beta\eta_{i+1,i},\, \alpha,\beta
> 0\}, \,i=1, \ldots, n-1.
\end{align*}

\begin{theorem}\label{th1}
Если рациональная функция \eqref{f1} удовлетворяет условиям
(\textbf{I}) --- (\textbf{III}), причем $P(z,w)\neq 0$ в вершинах
многоугольника $M$, не лежащих на координатных осях, тогда для
почти всех точек $(p,q)\in \mathbb{Z}^{2}_{+}$ на диагонали
$x=kp,\, y=kq,\, k\in\mathbb{Z}_{+}$ справедлива при $k \to
+\infty$ асимптотическая формула вида
\begin{align}
f(x,y) \sim
C(p,q;k)\frac{P(\hat{z},\hat{w})}{\hat{z}^{x+1}\hat{w}^{y+1}},
\end{align}
где для $(p,q)\in \Int{K_i}$ имеем $(\hat{z},\hat{w})=(z_i,w_i)$ и
$C(p,q;k)=\frac{const(p,q)}{\sqrt{k}}$, для $(p,q)\in \Int
\Omega_{ij}$ имеем $(\hat{z},\hat{w})=(z_{ij},w_{ij})$, а
$C(p,q;k)$ ---  не зависящая от $p,q,k$ константа.
\end{theorem}


Приведем пример того, что отказаться от условия \textbf{(III)} в
теореме \ref{th1} нельзя. Для функции
$F(z,w)=\frac{1}{(1-z-w)(1+z-w)}$ коэффициенты ее разложения в ряд
Тейлора имеют вид $f(x,y)=\frac{1}{2}
\binom{x+y+1}{x+1}\left[1+(-1)^{x}\right].$ Для любой точки
$(p,q)\in\mathbb{Z}^2_+$ найдется бесконечное число значений $k$
таких, что $f(kp,kq)=0$ (при любом $k=2,4,6,\ldots$),
следовательно, теорема \ref{th1} неверна.

В работе \cite{LPZ} введено понятие вектора Горна, который для
двойной последовательности $f(x,y)$ представляется в виде
$$\left( \frac{f(x+1,y)}{f(x,y)},\frac{f(x,y+1)}{f(x,y)} \right).$$

\begin{corollary} \label{HORN}
В условиях теоремы \ref{th1} почти для любого $(p,q)$ имеет место
равенство
\begin{align*}
    \lim_{\substack{x=kp\\y=kq\\k\to \infty}} \left( \frac{f(x+1,y)}{f(x,y)},\frac{f(x,y+1)}{f(x,y)}
    \right)=\left(\frac{1}{\hat{z}},\frac{1}{\hat{w}}\right),
\end{align*}
где $(\hat{z},\hat{w})\in\mathcal{V}$.
\end{corollary}

Доказательство теоремы \ref{th1} разобьем на ряд вспомогательных
утверждений.

\begin{lemma} \label{lemma1}
Пусть знаменатель рациональной функции \eqref{f1} имеет вид
$Q=Q_1Q_2$ и $P(z,w)\equiv 1$, тогда асимптотика коэффициентов
Тейлора $f(x,y)$ рациональной функции $F(z,w)$ на
$(p,q)$-диагонали $x=kp, y=kq$ при $k\to +\infty$ имеет вид
\begin{align} \label{f:lemma1}
    f(x,y) \sim
    \begin{cases}
        \frac{c_1}{\sqrt{k}}\frac{1}{{z}_1^{x+1}{w}_1^{y+1}},
                    &\text{если $(p,q) \in \Int K_1$},\\
        \frac{-1}{\Delta_{12}} \frac{1}{z_{12}^{x+1}w_{12}^{y+1}},
                    &\text{если $(p,q) \in \Int \Omega_{12}$,}\\
        \frac{c_2}{\sqrt{k}}\frac{1}{{z}_2^{x+1}{w}_2^{y+1}},
                    &\text{если $(p,q) \in \Int K_2$,}
    \end{cases}
\end{align}
где $c_i=\frac{(-1)^i}{\sqrt{2\pi}}
\frac{q}{pb_i(a_1-a_2)+qa_i(b_1-b_2)}\sqrt{\frac{p}{p+q}}, i=1,2$.
\end{lemma}

\textbf{Доказательство леммы \ref{lemma1}.} Воспользуемся формулой
Коши для коэффициентов Тейлора функции $F(z,w)$ и получим
\[f(x,y)=\int\limits_{\Gamma} \frac{1}{(1-a_1z-b_1w)(1-a_2z-b_2w)}
\frac{dz \wedge dw}{z^{x+1} w^{y+1}}, \text{где }
\Gamma=\{|z_1|=\varepsilon_1, |z_2|=\varepsilon_2\}.\]
Проинтегрировав по переменной $w$, получим
\begin{align*}
f(x,y)= -\frac{1}{\Delta_{12}}\int\limits_{|z|=\varepsilon}
\frac{1}{z-z_{12}} \left(\frac{b_1}{1-a_1z}\right)^{y+1}
\frac{dz}{z^{x+1}} + \frac{1}{\Delta_{12}}
\int\limits_{|z|=\varepsilon} \frac{1}{z-z_{12}}
\left(\frac{b_2}{1-a_2z}\right)^{y+1} \frac{dz}{z^{x+1}},
\end{align*}
где $\Delta_{12}$ --- определитель системы линейных уравнений
$Q_1=Q_2=0$. Обозначим
\begin{align} \label{phi_i}
\varphi_i(x,y)=\frac{1}{\Delta_{12}} \int\limits_{|z|=\varepsilon}
\frac{1}{z-z_{12}} \left(\frac{b_i}{1-a_i z}\right)^{y+1}
\frac{dz}{z^{x+1}},
\end{align}
тогда $f(x,y)=-\varphi_1(x,y)+\varphi_2(x,y)$.

Если $z_{ij}>0, w_{ij}>0$, то асимптотику интеграла \eqref{phi_i}
на $(p,q)$-диагонали $x=kp, y=kq$ при $k\to +\infty$ можно
вычислить, используя метод перевала (см., например, \cite{fedo}):
\begin{align*}
  \varphi_i(x,y) \sim
    \begin{cases}
        \hat{\varphi}_i(x,y) = \frac{c_i(p,q)}{\sqrt{k}} \frac{1}{z_i^{x+1}w_i^{y+1}},
            &\text{если $\frac{p}{q}<\frac{a_iz_{12}}{b_iw_{12}}$,}\\
        \frac{1}{\Delta_{12}}\frac{1}{z_{12}^{x+1} w_{12}^{y+1}},
            &\text{если $\frac{p}{q}>\frac{a_iz_{12}}{b_iw_{12}}$},
    \end{cases}
\end{align*}
где $z_i= \frac{1}{a_i}\frac{p}{p+q},
w_i=\frac{1}{b_i}\frac{q}{p+q}, i=1,2$ и $c_i(p,q)$ --- некоторая
константа.

Конусы \eqref{cones} в условиях леммы \ref{lemma1} можно
представить в виде $K_1=\{(p,q)\in \mathbb{Z}_+^2: 0\leqslant
\frac{p}{q}\leqslant \frac{a_1z_{12}}{b_1w_{12}}\},
\Omega_{12}=\{(p,q)\in \mathbb{Z}_+^2:
\frac{a_1z_{12}}{b_1w_{12}}\leqslant \frac{p}{q}\leqslant
\frac{a_2z_{12}}{b_2w_{12}}\}, K_2=\{(p,q)\in \mathbb{Z}_+^2:
\frac{a_2z_{12}}{b_2w_{12}}\leqslant \frac{p}{q}< +\infty\}$. Не
теряя общности, можно считать, что выполняются неравенства
$\frac{a_1}{b_1}<\frac{a_2}{b_2}$ и
$\frac{a_1(b_2-b_1)}{b_1(a_1-a_2)}<\frac{\ln b_2-\ln b_1}{\ln a_1
-\ln a_2}<\frac{a_2(b_2-b_1)}{b_2(a_1-a_2)}$, поэтому конус
$\Omega_{12}$ разбивается лучом $\{(p,q)\in\mathbb{R}^2:
\frac{p}{q}=\frac{\ln b_1 - \ln b_2}{\ln a_2 - \ln a_1}\}$ на
конусы $\Omega_{12}^1$ и $\Omega_{12}^2$, при этом если $(p,q)\in
\Int (K_1\cup\Omega_{12}^1)$, то
$\hat{\varphi}_2=o(\hat{\varphi}_1), k\to+\infty$ и если $(p,q)\in
\Int(\Omega_{12}^2\cup{K_2})$, то
$\hat{\varphi}_1=o(\hat{\varphi}_2), k\to+\infty$. Далее легко
показать, что для $(p,q)\in \Int K_1$ имеет место соотношение $f
\sim -\widehat{\varphi}_1$, так как $\varphi_1 \sim
\widehat{\varphi}_1, \varphi_2 \sim \widehat{\varphi}_2$ и
$\frac{\widehat{\varphi}_2}{\widehat{\varphi}_1}\to 0$ при $k \to
\infty$. Поэтому $\frac{f}{\widehat{\varphi}_1} =
-\frac{\varphi_1}{\widehat{\varphi}_1} +
\frac{\varphi_2}{\widehat{\varphi}_1} = -1+
\frac{\widehat{\varphi}_2}{\widehat{\varphi}_1} \to -1$.
Аналогично можно показать, что для  $(p,q)\in \Int K_2$  и
$(p,q)\in \Int\Omega_{12}$ справедливо $f \sim
\widehat{\varphi}_2$ и $f\sim -
\frac{1}{\Delta_{12}}\frac{1}{z_{12}^{x+1} w_{12}^{y+1}}$
соответственно.

Отметим, что если $z_{12}<0$, то по определению конус $K_2$
совпадет с $\mathbb{Z}^2_+$, а остальные конусы будут пустыми
множествами. Аналогично, в случае $w_{12}<0$, конус $K_1$ совпадет
с $\mathbb{Z}^2_+$, а остальные конусы будут пустыми множествами.
Лемма доказана.

\textbf{Пример.} Рассмотрим задачу о подбрасывании несимметричной
монеты. Вероятности выпадения аверса и реверса для первых $N$
испытаний равны $\frac{2}{3}$ и $\frac{1}{3}$ соответственно, а
для всех последующих испытаний --- $\frac{1}{3}$ и $\frac{2}{3}$
соответственно. Игрок желает получить $r$ аверсов и $s$ реверсов и
может выбрать $N$. В среднем сколько выборов $N\leqslant r+s$
будут желаемыми?

Вероятность, что $N$ будет желаемым выбором для игрока, есть сумма
по $a+b=N$ вида
\[
\binom{N}{a}\left(\frac{2}{3}\right)^a \left(\frac{1}{3}\right)^b
\binom{r+s-N}{r-a} \left(\frac{1}{3}\right)^{r-a}
\left(\frac{2}{3}\right)^{s-b}.
\]
Пусть $a_{rs}$ будет этим выражением, суммированным по $N$.
Последовательность $\{a_{rs}\}_{r,s\geqslant 0}$ является сверткой
двух последовательностей $\binom{r+s}{r}\left(\frac{2}{3}\right)^r
\left(\frac{1}{3}\right)^s$ и
$\binom{r+s}{r}\left(\frac{1}{3}\right)^r
\left(\frac{2}{3}\right)^s$. Производящая функция $F(z,w)$ этой
последовательности имеет вид
\[F(z,w)=\sum\limits_{r,s\geqslant 0} a_{rs}z^rw^s =
\frac{1}{(1-\frac{1}{3}z-\frac{2}{3}w)(1-\frac{2}{3}z-\frac{1}{3}w)}.\]
По лемме \ref{lemma1} получим $z_{12}=w_{12}=1$,
$K_1=\{(p,q):0\leqslant \frac{p}{q}\leqslant \frac{1}{2}\}$,
$\Omega_{12}=\{(p,q):\frac{1}{2}\leqslant \frac{p}{q}\leqslant
2\}$, $K_2=\{(p,q):2\leqslant \frac{p}{q}<+\infty\},
z_1=\frac{3p}{p+q}, w_1=\frac{3q}{2(p+q)}, z_2=\frac{3p}{2(p+q)},
w_2=\frac{3q}{p+q},
c_1=\frac{1}{\sqrt{2\pi}}\sqrt{\frac{p}{p+q}}\frac{9q}{q-2p},
c_2=\frac{1}{\sqrt{2\pi}}\sqrt{\frac{p}{p+q}}\frac{9q}{2q-p},
\Delta=a_1b_2-a_2b_1$, откуда при $k\to+\infty$
\begin{align*}
    f(x,y) \sim
    \begin{cases}
        \frac{1}{\sqrt{2\pi k}}\sqrt{\frac{p}{p+q}}\frac{2(p+q)^2}{(2p-q)p}
                    \left[\left(\frac{p+q}{3p}\right)^p\left(\frac{2(p+q)}{3q}\right)^q\right]^k,
                    &\text{если $(p,q) \in \Int K_1$},\\
        \frac{1}{a_2b_1-a_1b_2},
                    &\text{если $(p,q) \in \Int \Omega_{12}$,}\\
        \frac{1}{\sqrt{2\pi k}}\sqrt{\frac{p}{p+q}}\frac{2(p+q)^2}{(p-2q)p}
                    \left[\left(\frac{2(p+q)}{3p}\right)^p\left(\frac{p+q}{3q}\right)^q\right]^k,
                    &\text{если $(p,q) \in \Int K_2$.}
    \end{cases}
\end{align*}

\begin{lemma} \label{lemma2} Пусть задана рациональная функция двух
переменных
\begin{align} \label{f3}
F(z,w)=\frac{1}{\prod_{i=1}^{m}Q_i(z,w)}
\end{align}
и система многочленов $\{Q_i\}_{i=1}^{m}$ находится в общем
положении, тогда справедливо равенство
\begin{align} \label{razlozhenie}
F(z,w)= \sum_{1\leqslant i<j\leqslant m} \frac{A_{ij}}{Q_iQ_j},
\end{align}
где $A_{ij}$ некоторые константы.
\end{lemma}
\textbf{Доказательство леммы \ref{lemma2}.} Сначала докажем лемму
для случая $m=3$. Так как система $\{Q_i\}_{i=1}^{3}$ находится в
общем положении, то найдутся такие числа $A_{12}, A_{23}, A_{31}$,
что имеет место тождество $A_{12}Q_3+A_{23}Q_1+A_{31}Q_2\equiv 1$.
Действительно, это тождество эквивалентно системе трех линейных
уравнений с определителем
\[\Delta_{123}=
\begin{vmatrix}
  1   &   1 & 1   \\
  a_1 & a_2 & a_3 \\
  b_1 & b_2 & b_3 \\
\end{vmatrix}
\]
Разделив обе части этого  тождества на $Q_1Q_2Q_3$, получим
\begin{align} \label{form1}
\frac{1}{Q_1Q_2Q_3}=
\frac{A_{12}}{Q_1Q_2}+\frac{A_{23}}{Q_2Q_3}+\frac{A_{13}}{Q_1Q_3}.
\end{align}

Чтобы доказать лемму для случая $m=4$, домножим обе части
равенства \eqref{form1} на $\frac{1}{Q_4}$:
\[
\frac{1}{Q_1Q_2Q_3Q_4}=
\frac{1}{\Delta_{123}}\left(\frac{\Delta_{12}}{Q_1Q_2Q_4}+\frac{\Delta_{23}}{Q_2Q_3Q_4}+
\frac{\Delta_{13}}{Q_1Q_3Q_4}\right),
\]
и к каждому из слагаемых, стоящих в правой части, применим данную
лемму, доказанную для $m=3$. Утверждение леммы в случае
произвольного $m$ получается по индукции.

\begin{notation}
Константы $A_{ij}$ в разложении \eqref{razlozhenie} рациональной
функции \eqref{f3} на простейшие дроби можно вычислить по формуле
\[
A_{ij}=\frac{1}{(2\pi i)^2}\int\limits_{\Gamma_{ij}}
\frac{dz\wedge dw}{\prod_{i=1}^{m}Q_i(z,w)},
\]
где цикл $\Gamma_{ij}=\{|z-z_{ij}|=\varepsilon,
|w-w_{ij}|=\varepsilon\}$ и $(z_{ij},w_{ij})$ - решение системы
$Q_i=Q_j=0$.
\end{notation}

Утверждение леммы \ref{lemma2} является следствием общей теоремы
(см. \cite{lein}, \cite{uzhakov}) о разложении рациональной
функции на простейшие дроби. Однако в данном случае с точки зрения
приложений представляет интерес не столько возможность такого
разложения, сколько алгоритм, позволяющий найти коэффициенты
$A_{ij}$.

Асимптотическое поведение коэффициентов Тейлора $f(x,y)$ функции
$F(z,w)$ выражает

\begin{theorem} \label{th3}
Для рациональной функции вида \eqref{f3} почти для любой
$(p,q)$-диагонали $x=kp, y=kq, k\to +\infty$ справедлива
асимптотическая формула
\begin{align} \label{f10}
f(x,y) \sim
    \begin{cases}
        \frac{c_i}{\sqrt{k}}\frac{1}{z_i^{x+1}w_i^{y+1}},
                    &\text{если $(p,q)\in K_i$}\\
        \frac{A_{i,i+1}}{\Delta_{i,i+1}} \frac{1}{z_{i,i+1}^{x+1}w_{i,i+1}^{y+1}},
                    &\text{если $(p,q)\in \Omega_{i,i+1}$}
    \end{cases},
\end{align}
где $z_i=\frac{1}{a_i}\frac{p}{p+q},
w_i=\frac{1}{b_i}\frac{q}{p+q}$ --- решение системы
$Q_i=qzQ_z-pwQ_w=0$, $\Delta_{i,i+1}$ - определитель системы
линейных уравнений $Q_i=Q_{i+1}=0$, $(z_{i,i+1},w_{i,i+1})$ ---
решение этой системы, а $c_i$ --- некоторая констан-та.
\end{theorem}

\textbf{Доказательство теоремы \ref{th3}.} По лемме \ref{lemma2}
для функции \eqref{f3} имеет место разложение \eqref{razlozhenie},
тогда ее коэффициенты можно представить в виде суммы
\[
f(x,y) = \sum_{1\leqslant i<j\leqslant m} f_{ij}(x,y) =
\frac{1}{(2\pi i)^2} \sum_{1\leqslant i<j\leqslant m}
\int\limits_\Gamma \frac{A_{ij}}{Q_iQ_j} \frac{dz \wedge
dw}{z^{x+1}w^{y+1}}.
\]
К каждому слагаемому $f_{ij}$ применяем лемму \ref{lemma1}, тогда
почти для любой точки $(p,q)\in \mathbb{Z}^2_+$ либо
$f_{ij}(x,y)\sim
\frac{c_{ij}^{(i)}}{\sqrt{k}}\frac{1}{z_i^{x+1}w_i^{y+1}}$, либо
$f_{ij}(x,y)\sim
\frac{A_{ij}}{\Delta_{ij}}\frac{1}{z_{ij}^{x+1}w_{ij}^{y+1}}$ в
зависимости от конуса, которому принадлежит $(p,q)$. Среди этих
выражений необходимо выбрать главный член асимптотического
разложения.

Покажем, что если $(p,q)\in \Int K_i$, то справедлива
асимптотическая формула
$\frac{1}{z_i^{kp}w_j^{kq}}=o\left(\frac{1}{z_j^{kp}w_j^{kq}}\right)$
при $j\neq i$. Условие $(p,q)\in \Int K_i$ равносильно неравенству
$\frac{a_i(b_i-b_{i-1})}{b_i(a_{i-1}-a_i)} <\frac{p}{q}
<\frac{a_i(b_{i+1}-b_i)}{b_i(a_i-a_{i+1})}$.  Если $j<i$, то из
неравенства $\frac{\ln \frac{b_j}{b_i}}{\ln
\frac{a_i}{a_j}}<\frac{a_i(b_i-b_{i-1})}{b_i(a_{i-1}-a_i)}$
следует, что $\frac{p}{q}>\frac{\ln \frac{b_j}{b_i}}{\ln
\frac{a_i}{a_j}}$, а если $j>i$, то из неравенства $\frac{\ln
\frac{b_j}{b_i}}{\ln
\frac{a_i}{a_j}}>\frac{a_i(b_{i+1}-b_i)}{b_i(a_i-a_{i+1})}$
получаем $\frac{p}{q}<\frac{\ln \frac{b_j}{b_i}}{\ln
\frac{a_i}{a_j}}$, откуда $p \ln \frac{a_i}{a_j}> q \ln
\frac{b_j}{b_i}$ или
$\left(\frac{a_i}{a_j}\right)^p>\left(\frac{b_j}{b_i}\right)^q$
для всех $i\neq j$, а это значит, что $z_i^pw_i^q<z_j^pw_j^q$,
отсюда следует соотношение
$\frac{1}{z_i^{kp}w_j^{kq}}=o\left(\frac{1}{z_j^{kp}w_j^{kq}}\right)$.
Также легко видеть, что для $(p,q)\in \Int  K_i$, то
$\frac{1}{z_i^{kp}w_i^{kq}}=o\left(\frac{1}{z_{ij}^{kp}w_{ij}^{kq}}\right)$
при $j\neq i$. Аналогично можно показать, что если $(p,q)\in \Int
 \Omega_{i,i+1}$, то
$\frac{1}{z_{i,i+1}^{kp}w_{i,i+1}^{kq}}=o\left(\frac{1}{z_j^{kp}w_j^{kq}}\right)$.

\textbf{Доказательство теоремы \ref{th1}.} Рассмотрим интегральное
представление коэффициентов Тейлора фун\-кции \eqref{f1},
отвечающей условиям (\textbf{I}) --- (\textbf{III}). Пусть
$P(z,w)=\sum_{\alpha\beta} d_{\alpha\beta}z^\alpha w^\beta$
--- многочлен, тогда
\begin{align*}
f(x,y)= \frac{1}{(2\pi i)^2} \int\limits_\Gamma
\frac{P(z,w)}{\prod_i{Q_i}} \frac{dz\wedge dw}{z^{x+1}w^{y+1}}
=\frac{1}{(2\pi i)^2} \int\limits_\Gamma \frac{\sum_{\alpha\beta}
d_{\alpha\beta} z^\alpha
w^\beta}{\prod_i{Q_i}}\frac{dz\wedge dw}{z^{x+1}w^{y+1}}=\\
=\frac{1}{(2\pi i)^2} \sum_{\alpha \beta} d_{\alpha\beta}
\int\limits_\Gamma \frac{z^\alpha
w^\beta}{\prod_i{Q_i}}\frac{dz\wedge dw}{z^{x+1}w^{y+1}} =
\sum_{\alpha\beta} d_{\alpha\beta} g(x-\alpha, y-\beta).
\end{align*}
Из теоремы \ref{th3} следует, что $g(x-\alpha, y-\beta) \sim
C(p,q;k)\frac{\hat{z}^\alpha
\hat{w}^\beta}{\hat{z}^{x+1}\hat{w}^{y+1}}$, где значения
($\hat{z},\hat{w})$ и $C(p,q;k)$ определяются в соответствии с
формулой \eqref{f10}. Откуда
\begin{align*}
f(x,y)\sim \sum d_{\alpha\beta}\hat{z}^\alpha\hat{w}^\beta
c(p,q;k) \frac{1}{\hat{z}^{x+1}\hat{w}^{y+1}}=c(p,q;k)
\frac{P(\hat{z},\hat{w})}{\hat{z}^{x+1}\hat{w}^{y+1}},
\end{align*}
если $P(\hat{z},\hat{w})\neq 0$, а это справедливо почти для всех
$(p,q)\in\mathbb{Z}^2_+$. Теорема доказана.


\begin{thebibliography}{9}
    \bibitem{tz91} {\scshape Цих А.К.} {\itshape Условия абсолютной сходимоти ряда из коэффициентов
    Тейлора мероморфных функций двух переменных} / А.К. Цих // Мат. сб. --
    1991. -- № 11. -- С. 1588-1612.

    \bibitem{orlov2} {\scshape Орлов А.Г.} {\itshape Об асимптотике коэффициентов Тейлора рациональных
    функций двух переменных} / А.Г. Орлов // Изв. вузов. -- 1993. -- № 6. -- С. 26-33.

    \bibitem{orlovN} {\scshape Орлов А.Г.} {\itshape Об асимптотике коэффициентов Тейлора рациональных
    функций многих переменных} / А.Г. Орлов // Многомерный комплексный
    анализ. Межвузовский сборник. -- Красноярск. -- 1994. -- С. 116-141.

    \bibitem{makosi} {\scshape Макосий А.И.} {\itshape К вопросу об асимптотике коэффициентов
    ряда Тейлора} / А.И. Макосий // Многомерный комплексный анализ. -- Красноярск. --
    ИФ СО АН СССР. -- 1985. -- С. 224-227.

    \bibitem{pem}{\scshape Pemantle R., Wilson M.}  {\itshape Asymptotics for multivariate
    sequences, part I: smooth points of the singular variety} / R. Pemantle, M. Wilson // J.
    Comb. Th. -- Series A97. -- 129-161.

    \bibitem{LPZ} {\scshape Лейнартас Е.К., Пассаре М., Цих А.К.} {\itshape Асимптотика
    многомерных разностных уравнений} / Е.К. Лейнартас, М. Пассаре, А.К. Цих // УМН. --
    2005. -- Т. 60. -- Вып. 5(365). -- С. 171-172.

    \bibitem{fedo} {\scshape Федорюк М.В.} {\itshape Асимптотика. Интегралы и
    ряды} / М.В. Федорюк. -- М.: Наука, 1987. -- С. 295-297.

    \bibitem{lein} {\scshape Лейнартас Е.К.} {\itshape О разложении рациональных
    функций многих переменных на простейшие дроби} / Е.К. Лейнартас // Изв.
    вузов. -- Математика. -- 1978. -- № 10. -- С. 47-52.

    \bibitem{uzhakov} {\scshape Южаков А.П.} {\itshape Достаочное условие разделения аналитических
    особенностей в $\mathbb{C}^n$ и базис одного пространства голоморфных
    функций} / А.П. Южаков //  Матем. заметки. -- 1972. -- Т. 11. -- №~5. -- С. 585-596.

\end{thebibliography}
\end{document}